\documentclass[largeformat]{interact}

\usepackage[T1]{fontenc}
\usepackage[utf8]{inputenc}
\usepackage[english]{babel}
\usepackage{epstopdf}
\usepackage{mathtools,bm,mathrsfs}
\usepackage{array,enumitem,float}
\usepackage{natbib}
\bibpunct[, ]{(}{)}{;}{a}{}{,}

\usepackage[noperiod]{jabbrv}

\usepackage{placeins}

\renewcommand{\JournalTitle}[1]{#1}
\usepackage{url}
\usepackage[colorlinks=true,linkcolor=blue,citecolor=blue,urlcolor=blue]{hyperref}
\usepackage[nameinlink,capitalize]{cleveref}

\theoremstyle{plain}
\newtheorem{theorem}{Theorem}[section]
\newtheorem{proposition}[theorem]{Proposition}

\theoremstyle{definition}
\newtheorem{definition}{Definition}
\newtheorem{example}{Example}

\theoremstyle{remark}
\newtheorem{remark}{Remark}

\newcommand{\R}{\mathbb{R}}

\newcommand{\D}{\mathbb{D}}
\newcommand{\Dpos}{\mathbb{D}^{+}}

\newcommand{\kh}{\mathbf{k}}

\newcommand{\eone}{\mathbf{e}_{1}}
\newcommand{\etwo}{\mathbf{e}_{2}}

\newcommand{\preceqD}{\preceq_{\mathbb{D}}}

\newcommand{\Dint}[2]{[#1,#2]_{\D}}

\newcommand{\Ddist}{d_{h}}

\newcommand{\Edist}{d_{E}}

\newcommand{\PD}{\mathcal{P}_{\mathbb{D}}}

\newcommand{\maj}{\succcurlyeq_{\mathrm{maj}}}
\newcommand{\Dmaj}{\succcurlyeq_{\mathbb{D},\mathrm{maj}}}

\newcommand{\Neg}{\mathcal{N}}

\newcommand{\Yager}{\mathcal{Y}}

\title{Fixed-Point Results and Disequilibrium for Negators of Hyperbolic-Valued Probability Distributions}

\author{	\name{Edil D. Molina\textsuperscript{a},
		Juan Bory-Reyes\textsuperscript{b} and
		Jos\'e M. Sigarreta\textsuperscript{a}		\thanks{CONTACT Jos\'e M. Sigarreta.
			Email: 14366@uagro.mx}}	\affil{		\textsuperscript{a}Facultad de Matem\'aticas,
		Universidad Aut\'onoma de Guerrero,
		Carlos E. Adame No. 54, Col. Garita,
		Acapulco, Guerrero 39650, Mexico;\\
		\textsuperscript{b}Secci\'on de Estudios de
		Posgrado e Investigaci\'on,
		Escuela Superior de Ingenier\'ia Mec\'anica
		y El\'ectrica, Instituto Polit\'ecnico Nacional,
		Unidad Zacatenco, Mexico City, Mexico	}}

\begin{document}
\maketitle

\begin{abstract}
		Hyperbolic-valued probability distributions provide a probabilistic setting in which mass, order, and uncertainty take values in the positive cone of the hyperbolic algebra. We study the fixed-point structure of hyperbolic probability negators and the behavior of hyperbolic-valued disequilibrium along their iterates. Starting from a hyperbolic-valued metric on each mass stratum, we define the associated disequilibrium as the separation between a distribution and the uniform distribution determined by its hyperbolic mass. We consider the disequilibria induced by the hyperbolic Euclidean, \(\ell_1\), and Jensen--Shannon metrics. Their hyperbolic Schur-convexity is established and sharp lower and upper bounds are obtained. We also characterize the possible fixed-point sets of hyperbolic negators through an order-theoretic condition and derive criteria ensuring uniqueness of the associated uniform distribution. The behavior of disequilibrium along the iterates of hyperbolic negators is then analyzed through two complementary mechanisms. Domination by Yager's negator yields exact decay factors for the Euclidean and \(\ell_1\) disequilibria and a uniform decay factor for the Jensen--Shannon disequilibrium. Independently, Lipschitz conditions for generated negators lead to positive linear comparison maps whose spectral radii control the asymptotic decay of disequilibrium. These results provide structural and quantitative conditions ensuring that hyperbolic-valued disequilibrium converges to zero under repeated probability negation.
	\end{abstract}


\section{Introduction}
	\label{sec:introduction}
	
	Hyperbolic numbers form a commutative algebra endowed with a natural
	positive cone and the corresponding partial order. Their algebraic,
	analytic, and geometric properties have been studied in hyperbolic
	function theory, hypercomplex analysis, and related settings
	\citep{Cockle1848,Vignaux1935,Fjelstad1986,Sobczyk1995,Kisil2012,Sobczyk2013}.
	In probability theory, the hyperbolic setting has led to
	hyperbolic-valued measures, entropy, and other probabilistic
	constructions
	\citep{Alpay2017,TellezSanchez2021,Alpay2025}. Hyperbolic-valued
	seminorms, metrics, and topological structures have also been developed
	in bicomplex and hyperbolic analysis 
	\citep{LunaElizarraras2014,Kumar2016,Ghosh2023,AliKumar2026,Saini2020}. These
	developments provide tools for studying probabilistic and dynamical
	phenomena while preserving the ordered hyperbolic structure.
	
	The negation of probability distributions was introduced by
	\citet{Yager2015} as a mechanism for transforming a probability
	distribution toward a state of greater uncertainty. Subsequent developments have
	addressed structural properties, fixed points, contraction,
	involutivity, generated and parametric negators, and information loss
	\citep{Srivastava2017,Srivastava2021,Kaur2022,Klein2022,
		Batyrshin2021a,Batyrshin2021b,Batyrshin2021,Batyrshin2023,
		EnsasteguiOrtega2024}. Generalized and information-theoretic forms of
	probability negation, together with applications in information fusion,
	recognition, decision making, and evidence theory, have also been widely
	studied
	\citep{Srivastava2019,Zhang2019,Liu2021,Gao2019,Sun2020,Pham2021,
		Wu2022,Tanwar2023,Tanwar2023a,Yin2019,Luo2020,Deng2020,Liu2023,BrunoVaccaro2026}.
	
	The theory of negation for finite hyperbolic-valued probability
	distributions was initiated by \citet{BoryReyes2026}, who introduced
	foundational notions such as hyperbolic negators, generated negators,
	Yager's hyperbolic negator, and a hyperbolic majorization framework. For a broad class of
	generated negators, majorization was used to establish entropy increase
	and quantitative convergence under iteration. The hyperbolic setting
	also revealed a phenomenon absent from the classical theory, since a
	negator may possess fixed points distinct from the associated uniform
	distribution.
	
	The notion of disequilibrium has a well-established role in the study of
	probability distributions. In the statistical-complexity framework
	introduced by 
	\citet{LopezRuiz1995}, disequilibrium measures the departure
	of a probability distribution from an equilibrium distribution. This viewpoint was subsequently
	developed using different notions of distance in probability space,
	including constructions based on the Jensen--Shannon divergence
	\citep{Martin2003,Lamberti2004,Kowalski2011}. Motivated by this approach, we consider
	disequilibrium in the hyperbolic-valued setting. Given a
	hyperbolic-valued metric \(d\) on each mass stratum of hyperbolic-valued probability distribution \(P\in \PD(n)\), the
	associated disequilibrium is defined by
	\(D_d(P)=d(P,U_n^{s(P)})\), where \(U_n^{s(P)}\) is determined by the
	hyperbolic mass of \(P\). The metric supplies the underlying notion of
	separation, while the resulting map \(D_d\) provides the quantity used
	throughout the subsequent analysis.
	
	We consider three instances of this construction. The hyperbolic
	Euclidean and \(\ell_1\) metrics induce geometric disequilibria, while
	the hyperbolic Jensen--Shannon metric induces an information-theoretic
	one. The latter is obtained from the Jensen--Shannon divergence and,
	ultimately, from the Kullback--Leibler divergence
	\citep{Kullback1951,Lin1991,Endres2003}. We prove that the resulting maps
	\(D_E\), \(D_h\), and \(D_{\operatorname{JS}}\) are hyperbolically
	Schur-convex, with strict Schur-convexity for \(D_E\) and
	\(D_{\operatorname{JS}}\). The extremal structure induced by
	hyperbolic majorization then yields sharp lower and upper bounds for
	each disequilibrium.
	
	The fixed-point problem is studied independently of these particular
	choices of disequilibrium. We introduce the class of order-antichain
	distributions and characterize exactly which subsets of a mass stratum
	can arise as the fixed-point set of a hyperbolic negator. This
	characterization shows that the associated uniform distribution is the
	only possible fixed point in the strata of masses \(\eone\) and
	\(\etwo\), whereas additional fixed points may occur in the stratum of
	mass \(1\). We then derive two types of sufficient conditions for
	uniqueness. The first is based on level sets of admissible hyperbolic
	generators and the second uses comparison maps acting on
	hyperbolic-valued disequilibrium.
	
	We next study the behavior of disequilibrium under repeated negation.
	A general comparison principle shows that an estimate of the form
	\(D_d(\Neg(P))\preceqD\varphi(D_d(P))\), together with convergence to
	zero of the iterates of \(\varphi\), implies
	\(D_d(\Neg^m(P))\to0\). Two mechanisms for obtaining such estimates are
	developed. The first uses hyperbolic majorization and domination by
	Yager's negator. For the Euclidean and \(\ell_1\) disequilibria,
	Yager's negator produces exact decay factors. For the Jensen--Shannon
	disequilibrium, strict Schur-convexity combined with a local analysis
	provides a uniform decay factor strictly smaller than one when \(n>2\).
	It follows that the three disequilibria converge to zero for every
	negator satisfying the corresponding Yager-domination condition.
	
	The second mechanism is independent of majorization and applies to
	generated negators satisfying suitable Lipschitz conditions. Lower
	bounds for the normalizing factors and Lipschitz constants of the
	generator produce positive real-linear comparison maps controlling the
	three disequilibria. These maps are represented by non-negative
	matrices, and conditions on their spectral radii ensure decay to zero.
	This provides a complementary connection between regularity properties
	of the generator and the asymptotic behavior of hyperbolic-valued
	disequilibrium.
	
	The paper is organized as follows. Section~\ref{sec:preliminaries}
	recalls the basic notions of hyperbolic-valued probability
	distributions, hyperbolic majorization, and probability negation.
	Section~\ref{sec:hyperbolic-valued-distances} introduces
	hyperbolic-valued metrics as the basis for defining disequilibrium and
	establishes the Schur-convexity and sharp bounds of
	\(D_E\), \(D_h\), and \(D_{\operatorname{JS}}\).
	Section~\ref{sec:uniform-fixed-point-uniqueness} studies the
	fixed-point structure of hyperbolic negators and sufficient conditions
	for uniqueness of the associated uniform distribution.
	Section~\ref{sec:disequilibrium-decay} develops general criteria for
	the decay of disequilibrium together with the Yager-domination and
	Lipschitz--spectral mechanisms.

	\section{Preliminaries}
	\label{sec:preliminaries}
	
	This section recalls the notation and basic results used throughout the
	paper within the finite hyperbolic-valued probability framework
	developed by \citet{BoryReyes2026}.
	
	The set of hyperbolic numbers is the commutative real algebra
	\[
	\D=\{a+\kh b\mid a,b\in\R,\ \kh^2=1\}.
	\]
	Its standard idempotent basis is
	\[
	\eone=\frac{1+\kh}{2},\qquad \etwo=\frac{1-\kh}{2},\qquad
	\eone^2=\eone,\quad \etwo^2=\etwo,\quad \eone\etwo=0,\quad \eone+\etwo=1.
	\]
	Every \(z\in\D\) has a unique representation \(z=z^{(1)}\eone+z^{(2)}\etwo\); if \(z=a+\kh b\), then \(z^{(1)}=a+b\) and \(z^{(2)}=a-b\). Addition and multiplication are componentwise in this basis.
	
	The natural partial order on \(\D\) is
	\[
	z\preceqD w\quad\Longleftrightarrow\quad
	z^{(1)}\le w^{(1)}\ \text{and}\ z^{(2)}\le w^{(2)}.
	\]
	We write
	\[
	\Dpos=\{z\in\D\mid 0\preceqD z\},\qquad
	\Dint{0}{1}
	=
	\{z^{(1)}\eone+z^{(2)}\etwo\mid 0\le z^{(1)},z^{(2)}\le1\}.
	\]
	The order \(\preceqD\) is not total; for instance, \(\eone\) and \(\etwo\) are incomparable. This feature will be relevant in the study of fixed points.
	
	Let \(n\ge2\), let \(\Omega=\{\omega_1,\ldots,\omega_n\}\), and write
	\[
	P=(p_1,\ldots,p_n)\in\D^n,\qquad
	p_j=p_j^{(1)}\eone+p_j^{(2)}\etwo,\qquad
	P^{(\ell)}=(p_1^{(\ell)},\ldots,p_n^{(\ell)}),\quad \ell=1,2.
	\]
	
	\begin{definition}
		\label{def:hyperbolic-probability-distribution}
		An \(n\)-tuple \(P=(p_1,\ldots,p_n)\in\D^n\) is called a hyperbolic-valued probability distribution on \(\Omega\) if \(p_j\in\Dint{0}{1}\) for every \(j\) and
		$
		\sum_{j=1}^{n}p_j\in\{1,\eone,\etwo\}
		$.
		
		The set of all such distributions is denoted by \(\PD(n)\). For \(P\in\PD(n)\), its total mass is
		\[
		s(P)=\sum_{j=1}^{n}p_j
		=s^{(1)}(P)\eone+s^{(2)}(P)\etwo,\qquad
		s^{(\ell)}(P)=\sum_{j=1}^{n}p_j^{(\ell)}.
		\]
	\end{definition}
	
	Since \(s^{(\ell)}(P)\in\{0,1\}\), the component \(P^{(\ell)}\) is either a real probability distribution or the zero vector; these components are called active and inactive, respectively. For \(s\in\{1,\eone,\etwo\}\), the uniform hyperbolic-valued probability distribution of length \(n\) and total mass \(s\) is
	\[
	U_n^s=\left(\frac{s}{n},\ldots,\frac{s}{n}\right),
	\]
	and the uniform distribution associated with \(P\in\PD(n)\) is \(U_n^{s(P)}\). 
	
	We call \(P\in\PD(n)\) a hyperbolic point distribution if each active idempotent component is a real point distribution.
	
	To compare the dispersion of hyperbolic-valued distributions, we use majorization componentwise. 
	
	Set
	\[
	\Delta_n=\left\{x=(x_1,\ldots,x_n)\in\R_{\ge0}^n\mid \sum_{j=1}^{n}x_j=1\right\},
	\]
	and let \(x^\downarrow\) denote the decreasing rearrangement of \(x\). For non-negative vectors \(x,y\in\R^n\) with the same total mass, we write \(x\maj y\) if for $k=1,\ldots,n-1$
	\[
	\sum_{j=1}^{k}x_j^\downarrow\ge\sum_{j=1}^{k}y_j^\downarrow,\quad \text{and}
	\qquad
	\sum_{j=1}^{n}x_j=\sum_{j=1}^{n}y_j.
	\]
	
	We use the fact that symmetric convex functions are Schur-convex,
	and that \(F(x)=\sum_{j=1}^{n}\varphi(x_j)\) is strictly Schur-convex
	whenever \(\varphi\) is strictly convex; see \citet{Hardy1934}
	and \citet[Chapter~3]{Marshall2011}.
	
	\begin{definition}
		\label{def:hyperbolic-majorization}
		Let \(P,Q\in\PD(n)\) with \(s(P)=s(Q)\). We say that \(P\) hyperbolically majorizes \(Q\), and write \(P\Dmaj Q\), if
		\[
		P^{(1)}\maj Q^{(1)}\qquad\text{and}\qquad P^{(2)}\maj Q^{(2)}.
		\]
	\end{definition}
	
	The relation \(\Dmaj\) is a preorder, and mutual majorization occurs precisely when the corresponding idempotent components differ only by permutations. Its extremal elements on each mass stratum are described by the
	following result established by \citet{BoryReyes2026}.
	
	\begin{proposition}
		\label{prop:vertex-uniform-majorization}
		Let \(P\in\PD(n)\). If \(V\in\PD(n)\) is a hyperbolic point distribution with \(s(V)=s(P)\), then
		\[
		V\Dmaj P\Dmaj U_n^{s(P)}.
		\]
	\end{proposition}
	
	A transformations that reverse the hyperbolic order of the coordinates of a hyperbolic probability distribution, while preserving total mass is defined as follows.
	
	\begin{definition}
		\label{def:hyperbolic-negator}
		A map \(\Neg:\PD(n)\to\PD(n)\) is called a hyperbolic negator if, for every \(P=(p_1,\ldots,p_n)\in\PD(n)\), with \(\Neg(P)=Q=(q_1,\ldots,q_n)\), one has \(s(Q)=s(P)\) and
		\[
		p_i\preceqD p_j\quad\to\quad q_j\preceqD q_i.
		\]
	\end{definition}
	
	If \(p_i=p_j\), applying the negator property in both directions gives \(q_i=q_j\). Consequently, every hyperbolic negator fixes the uniform distribution:
	\begin{equation}
		\label{eq:uniform-fixed-point}
		\Neg(U_n^s)=U_n^s,\qquad s\in\{1,\eone,\etwo\}.
	\end{equation}
	
	A broad class of hyperbolic negators is obtained by normalizing order-reversing generators.
	
	\begin{definition}
		\label{def:admissible-hyperbolic-generator}
		A map \(\Phi:\Dint{0}{1}\to\Dpos\), with
		\(\Phi(z)=\Phi^{(1)}(z)\eone+\Phi^{(2)}(z)\etwo\), is called an admissible hyperbolic generator if:
		\begin{enumerate}[label=\textup{(\roman*)}]
			\item \(z\preceqD w\) implies \(\Phi(w)\preceqD\Phi(z)\);
			\item for every \(P=(p_1,\ldots,p_n)\in\PD(n)\) and every active component \(\ell\),
			\[
			Z_\ell^\Phi(P)=\sum_{j=1}^{n}\Phi^{(\ell)}(p_j)>0.
			\]
		\end{enumerate}
	\end{definition}
	
	\begin{definition}
		\label{def:generated-hyperbolic-transformation}
		Let \(\Phi:\Dint{0}{1}\to\Dpos\) be an admissible hyperbolic generator. The generated transformation \(\Neg_\Phi:\PD(n)\to\PD(n)\), \(\Neg_\Phi(P)=Q=(q_1,\ldots,q_n)\), is defined by
		\[
		q_j^{(\ell)}
		=
		\begin{cases}
			\frac{\Phi^{(\ell)}(p_j)}{Z_\ell^\Phi(P)}, & s^{(\ell)}(P)=1,\\[8pt]
			0, & s^{(\ell)}(P)=0,
		\end{cases}
		\qquad
		q_j=q_j^{(1)}\eone+q_j^{(2)}\etwo.
		\]
	\end{definition}
	
	The normalization in the preceding definition preserves the total mass and the order-reversing property.
	
	\begin{proposition}
		\label{prop:generated-transformation-is-negator}
		Every admissible hyperbolic generator \(\Phi\) defines a hyperbolic negator \(\Neg_\Phi\).
	\end{proposition}
	
	The affine generator \(\Phi_Y(z)=1-z\) gives Yager's hyperbolic negator
	\[
	\Yager(P)=\left(\frac{s(P)-p_1}{n-1},\ldots,\frac{s(P)-p_n}{n-1}\right).
	\]
	Its position in the hyperbolic majorization order will provide the comparison principle used later.
	
	\begin{proposition}
		\label{prop:yager-intermediate-majorization}
		For every \(P\in\PD(n)\), one has \(P\Dmaj\Yager(P)\). Moreover, if \(\Neg:\PD(n)\to\PD(n)\) satisfies \(\Yager(P)\Dmaj\Neg(P)\) for every \(P\in\PD(n)\), then \(P\Dmaj\Neg(P)\).
	\end{proposition}
	
	In the classical setting, \citet{BrunoVaccaro2026} established an
	analogous majorization relation for independent probability negators,
	showing that Yager's negation majorizes the distributions obtained
	from this class.
	
	For generated negators, the following criterion ensures domination by Yager's negator.
	
	\begin{theorem}
		\label{thm:generated-negator-yager-dominated}
		Let \(\Phi:\Dint{0}{1}\to\Dpos\) be an admissible hyperbolic generator. Suppose that, for every \(P=(p_1,\ldots,p_n)\in\PD(n)\), every active component \(\ell\), and all \(i,j\in\{1,\ldots,n\}\),
		\[
		p_j^{(\ell)}=1\to\Phi^{(\ell)}(p_j)=0,
		\]
		and
		\[
		p_i^{(\ell)}\le p_j^{(\ell)}<1
		\to
		\Phi^{(\ell)}(p_i)\ge\Phi^{(\ell)}(p_j)
		\quad\text{and}\quad
		\frac{\Phi^{(\ell)}(p_i)}{1-p_i^{(\ell)}}
		\le
		\frac{\Phi^{(\ell)}(p_j)}{1-p_j^{(\ell)}}.
		\]
		Then, for every \(P\in\PD(n)\),
		\[
		\Yager(P)\Dmaj\Neg_\Phi(P)\qquad\text{and}\qquad P\Dmaj\Neg_\Phi(P).
		\]
	\end{theorem}

	\section{Hyperbolic-valued metrics and disequilibrium}
	\label{sec:hyperbolic-valued-distances}
	
	The notion of disequilibrium is used in probability space to quantify
	the departure of a probability distribution from the uniform equilibrium
	distribution. In the classical setting, this notion was introduced within
	the statistical-complexity framework and subsequently developed through
	different distance-based formulations
	\citep{LopezRuiz1995,Martin2003,Lamberti2004,Kowalski2011}.
	This viewpoint extends naturally to \(\PD(n)\), where every
	\(P\in\PD(n)\) has an associated uniform distribution
	\(U_n^{s(P)}\). We therefore use hyperbolic-valued metrics to measure
	this departure from equilibrium and refer to the resulting quantity as
	hyperbolic-valued disequilibrium.
	
	\begin{definition}
		\label{def:hyperbolic-valued-metric}
		Let \(X\) be a nonempty set. A map
		\[
		d:X\times X\to\Dpos
		\]
		is called a hyperbolic-valued metric if, for all \(x,y,z\in X\),
		\begin{enumerate}[label=\textup{(\roman*)}]
			\item \(0\preceqD d(x,y)\), and \(d(x,y)=0\) if and only if \(x=y\);
			\item \(d(x,y)=d(y,x)\);
			\item \(d(x,z)\preceqD d(x,y)+d(y,z)\).
		\end{enumerate}
	\end{definition}
	
	\begin{definition}
		\label{def:hyperbolic-valued-disequilibrium}
		Let \(d\) be a hyperbolic-valued metric on each mass stratum of
		\(\PD(n)\). The hyperbolic-valued disequilibrium associated with
		\(d\) is the map \(D_d:\PD(n)\to\Dpos\)
		defined by
		\[
		D_d(P)
		=
		d\left(P,U_n^{s(P)}\right).
		\]
	\end{definition}
	
	We will now introduce three hyperbolic-valued metrics on \(\PD(n)\) and the corresponding disequilibria induced by them.
	
	\begin{example}[Euclidean disequilibrium]
		\label{ex:hyperbolic-euclidean-metric}
		For \(P,Q\in\PD(n)\), define
		\[
		\Edist(P,Q)
		=
		\left(
		\sum_{j=1}^{n}
		\left(p_j^{(1)}-q_j^{(1)}\right)^2
		\right)^{1/2}\eone
		+
		\left(
		\sum_{j=1}^{n}
		\left(p_j^{(2)}-q_j^{(2)}\right)^2
		\right)^{1/2}\etwo.
		\]
		The metric axioms follow directly from those of the classical
		Euclidean metric. The associated hyperbolic-valued disequilibrium is
		\[
		D_E(P)
		=
		\Edist\left(P,U_n^{s(P)}\right).
		\]
	\end{example}
	
	\begin{example}[\(\ell^1\) disequilibrium]
		\label{ex:hyperbolic-l1-metric}
		For \(P,Q\in\PD(n)\), define
		\[
		\Ddist(P,Q)
		=
		\left(
		\sum_{j=1}^{n}
		\left|p_j^{(1)}-q_j^{(1)}\right|
		\right)\eone
		+
		\left(
		\sum_{j=1}^{n}
		\left|p_j^{(2)}-q_j^{(2)}\right|
		\right)\etwo.
		\]
		The metric axioms follow from those of the classical
		\(\ell^1\) metric. The associated hyperbolic-valued disequilibrium is
		\[
		D_h(P)
		=
		\Ddist\left(P,U_n^{s(P)}\right).
		\]
	\end{example}
	
	\begin{example}[Jensen--Shannon disequilibrium]
		\label{ex:hyperbolic-js-metric}
		For \(P,Q\in\PD(n)\) with \(s(P)=s(Q)\), the
		hyperbolic-valued Kullback--Leibler divergence of \(P\)
		relative to \(Q\) is defined by
		\[
		\operatorname{KL}_{\D}(P\|Q)
		=
		K_1(P\|Q)\eone
		+
		K_2(P\|Q)\etwo,
		\]
		where, for \(\ell\in\{1,2\}\),
		\[
		K_\ell(P\|Q)
		=
		\begin{cases}
			\sum_{j=1}^{n}
			p_j^{(\ell)}
			\log\frac{p_j^{(\ell)}}{q_j^{(\ell)}},
			&
			s^{(\ell)}(P)=1,
			\\[2.5ex]
			0,
			&
			s^{(\ell)}(P)=0.
		\end{cases}
		\]
		We adopt the conventions that \(0\log\left(\frac{0}{q}\right) =0\) and \(p\log\left( \frac{p}{0}\right) =+\infty\), for \(p>0\).
		
		Set \(M(P,Q)=\frac{P+Q}{2}\). The hyperbolic-valued Jensen--Shannon divergence is
		\[
		\operatorname{JS}_{\D}(P,Q)
		=
		\frac12
		\operatorname{KL}_{\D}
		\left(P\middle\|M(P,Q)\right)
		+
		\frac12
		\operatorname{KL}_{\D}
		\left(Q\middle\|M(P,Q)\right).
		\]
		Since \(M(P,Q)\) contains the support of both \(P\) and \(Q\),
		\(\operatorname{JS}_{\D}(P,Q)\in\Dpos\). The corresponding
		hyperbolic-valued Jensen--Shannon metric is
		\[
		d_{\operatorname{JS}}^{\D}(P,Q)
		=
		\sqrt{\operatorname{JS}_{\D}(P,Q)}.
		\]
		The fact that the square root of the classical Jensen--Shannon
		divergence is a metric \citep{Endres2003} yields the metric properties
		of \(d_{\operatorname{JS}}^{\D}\) on each mass stratum of
		\(\PD(n)\).
		
		The associated Jensen--Shannon disequilibrium is defined by
		\[
		D_{\operatorname{JS}}(P)
		=
		d_{\operatorname{JS}}^{\D}
		\left(P,U_n^{s(P)}\right).
		\]
	\end{example}
		
	We next study the behavior of \(D_E\), \(D_h\), and
	\(D_{\operatorname{JS}}\) under hyperbolic majorization.
	
	\begin{definition}
		\label{def:hyperbolic-schur-convexity}
		A map \(\mathcal F:\PD(n)\to\Dpos\) is called hyperbolically
		Schur-convex if
		\[
		P\Dmaj Q
		\quad\to\quad
		\mathcal F(Q)\preceqD\mathcal F(P).
		\]
		It is strictly hyperbolically Schur-convex if equality implies that,
		for every \(\ell\in\{1,2\}\) such that \(s^{(\ell)}(P)=1\),
		\(P^{(\ell)}\) and \(Q^{(\ell)}\) differ only by a permutation.
		Hyperbolic Schur-concavity is defined by reversing the inequality.
	\end{definition}
	
	Monotonicity with respect to hyperbolic majorization is established in the following result.
	
	\begin{proposition}
		\label{prop:uniform-metrics-schur-convex}
		The maps \(D_E\), \(D_h\), and
		\(D_{\operatorname{JS}}\) are hyperbolically Schur-convex on
		\(\PD(n)\). Moreover, \(D_E\) and \(D_{\operatorname{JS}}\) are
		strictly hyperbolically Schur-convex.
	\end{proposition}
	
	\begin{proof}
		Let \(P,Q\in\PD(n)\) satisfy \(P\Dmaj Q\). Then, \(s(P)=s(Q)\).
		Fix \(\ell\in\{1,2\}\) such that \(s^{(\ell)}(P)=1\).
		
		For \(D_E\) we have,
		\[
		\left(\left[D_E(P)\right]^{(\ell)}\right)^2
		=
		\sum_{j=1}^{n}
		\left(p_j^{(\ell)}-\frac1n\right)^2
		=
		\sum_{j=1}^{n}\left(p_j^{(\ell)}\right)^2-\frac1n.
		\]
		Since \(t\mapsto t^2\) is strictly convex, the right-hand side is
		strictly Schur-convex. As the square-root function is strictly
		increasing we get
		\[
		\left[D_E(Q)\right]^{(\ell)}
		\le
		\left[D_E(P)\right]^{(\ell)},
		\]
		with equality only when \(P^{(\ell)}\) and \(Q^{(\ell)}\) differ
		by a permutation.
		
		For \(D_h\) it follows that,
		\[
		\left[D_h(P)\right]^{(\ell)}
		=
		\sum_{j=1}^{n}
		\left|p_j^{(\ell)}-\frac1n\right|.
		\]
		Since \(t\mapsto |t-1/n|\) is convex, the right-hand side is a
		symmetric separable Schur-convex function. Hence
		\[
		\left[D_h(Q)\right]^{(\ell)}
		\le
		\left[D_h(P)\right]^{(\ell)}.
		\]
		
		Finally,
		\[
		\left[
		\operatorname{JS}_{\D}
		\left(P,U_n^{s(P)}\right)
		\right]^{(\ell)}
		=
		\sum_{j=1}^{n}\psi_n\left(p_j^{(\ell)}\right),
		\]
		where
		\[
		\psi_n(t)
		=
		-\frac12\left(t+\frac1n\right)
		\log\left(\frac{t+1/n}{2}\right)
		+
		\frac12t\log t
		-
		\frac{1}{2n}\log n.
		\]
		For \(0<t\le1\),
		\[
		\psi_n''(t)
		=
		\frac{1}{2t}
		-
		\frac{1}{2(t+1/n)}
		=
		\frac{1}{2nt(t+1/n)}
		>0.
		\]
		Thus \(\psi_n\) is strictly convex on \([0,1]\), and therefore
		\[
		\left[D_{\operatorname{JS}}(Q)\right]^{(\ell)}
		\le
		\left[D_{\operatorname{JS}}(P)\right]^{(\ell)},
		\]
		with equality only when \(P^{(\ell)}\) and \(Q^{(\ell)}\) differ
		by a permutation.
		
		If \(s^{(\ell)}(P)=0\), the corresponding coefficients vanish.
		Therefore,
		\[
		D_E(Q)\preceqD D_E(P),\qquad
		D_h(Q)\preceqD D_h(P),\qquad
		D_{\operatorname{JS}}(Q)\preceqD D_{\operatorname{JS}}(P).
		\]
	\end{proof}
	
	The preceding proposition, together with
	Proposition~\ref{prop:vertex-uniform-majorization}, yields sharp bounds
	for \(D_E\), \(D_h\), and \(D_{\operatorname{JS}}\). Set
	\[
	C_n
	=
	\log 2+\frac12\log n
	-
	\frac{n+1}{2n}\log(n+1).
	\]
	
	\begin{theorem}
		\label{thm:uniform-metric-sharp-bounds}
		Let \(n\ge2\) and \(P\in\PD(n)\). Then
		\[
		0
		\preceqD
		D_E(P)
		\preceqD
		\sqrt{1-\frac1n}\,s(P),
		\]
		\[
		0
		\preceqD
		D_h(P)
		\preceqD
		2\left(1-\frac1n\right)s(P),
		\]
		and
		\[
		0
		\preceqD
		D_{\operatorname{JS}}(P)
		\preceqD
		\sqrt{C_n}\,s(P).
		\]
		In each case, the lower bound is attained if and only if
		\(P=U_n^{s(P)}\), whereas the upper bound is attained if and only if
		\(P\) is a hyperbolic point distribution.
	\end{theorem}
	
	\begin{proof}
		Let \(V\in\PD(n)\) be a hyperbolic point distribution such that
		\(s(V)=s(P)\). By
		Proposition~\ref{prop:vertex-uniform-majorization},
		\[
		V\Dmaj P\Dmaj U_n^{s(P)}.
		\]
		Hence, by
		Proposition~\ref{prop:uniform-metrics-schur-convex},
		\[
		0
		=
		D_E\left(U_n^{s(P)}\right)
		\preceqD
		D_E(P)
		\preceqD
		D_E(V),
		\]
		and analogously
		\[
		0
		\preceqD
		D_h(P)
		\preceqD
		D_h(V),
		\qquad
		0
		\preceqD
		D_{\operatorname{JS}}(P)
		\preceqD
		D_{\operatorname{JS}}(V).
		\]
		
		A direct calculation gives
		\[
		D_E(V)
		=
		\sqrt{1-\frac1n}\,s(P),
		\qquad
		D_h(V)
		=
		2\left(1-\frac1n\right)s(P),
		\qquad
		D_{\operatorname{JS}}(V)
		=
		\sqrt{C_n}\,s(P).
		\]
		
		Suppose that \(D_E(P)=0\) and so \(d_E(P,U_n^{s(P)})=0\), which occurs if and only if
		\(P=U_n^{s(P)}\), because
		\(d_E\) is a hyperbolic-valued metric. The same argument applies to \(D_h\) and
		\(D_{\operatorname{JS}}\). Therefore, in each case, the lower bound
		is attained if and only if \(P=U_n^{s(P)}\).
		
		For \(D_E\) and \(D_{\operatorname{JS}}\), strict hyperbolic
		Schur-convexity implies that equality in the upper bound holds
		exactly when \(P\) is a hyperbolic point distribution.
		
		For \(D_h\), equality in the upper bound requires, for every
		\(\ell\in\{1,2\}\) such that \(s^{(\ell)}(P)=1\),
		\[
		\sum_{j=1}^{n}
		\left|p_j^{(\ell)}-\frac1n\right|
		=
		2\left(1-\frac1n\right).
		\]
		This equality holds if and only if \(P^{(\ell)}\) is a point
		distribution. Hence equality in the upper bound holds exactly when
		\(P\) is a hyperbolic point distribution.
	\end{proof}
	
	\section{Fixed-point uniqueness for hyperbolic negators}
	\label{sec:uniform-fixed-point-uniqueness}
	
	The study of generalized metrics with values in ordered vector spaces
	and cones goes back to the early work of Kurepa and Kantorovich, and
	was later developed in connection with Banach-type fixed-point
	principles by Schröder, Perov, and subsequent authors
	\citep{Kurepa1934,Kantorovich1939,Schroder1956,Perov1964,Zabreiko1997}.
	
	Hyperbolic-valued metrics provide a natural setting for analogous
	fixed-point questions \citep{Ghosh2023,Sager2021,Sager2023}.
	
	For hyperbolic negators, however, the fixed-point structure is also
	strongly influenced by the partial order of \(\D\). As observed by
	\citet{BoryReyes2026}, the possibility of fixed points other than
	the associated uniform distributions is a distinctive feature of the
	hyperbolic-valued setting. We examine this phenomenon
	by first describing the order-theoretic structure of fixed points and
	then deriving structural and metric conditions that ensure uniqueness
	of the associated uniform distribution.
	
	For a hyperbolic negator \(\Neg:\PD(n)\to\PD(n)\) and
	\(s\in\{1,\eone,\etwo\}\), we write
	\[
	\operatorname{Fix}_s(\Neg)
	=
	\{P\in\PD(n)\mid s(P)=s,\ \Neg(P)=P\}.
	\]
	By \eqref{eq:uniform-fixed-point}, 	\(U_n^s\in\operatorname{Fix}_s(\Neg)\).
	
	\medskip
	
	We begin with the order-theoretic restriction satisfied by every fixed
	point.
	
	\begin{definition}
		\label{def:order-antichain-distribution}
		Let \(P=(p_1,\ldots,p_n)\in\PD(n)\). We say that \(P\) is an
		order-antichain distribution if \(p_i\ne p_j\) implies
		\(p_i\not\preceqD p_j\) and \(p_j\not\preceqD p_i\). We denote by
		\[
		\mathcal A_n
		=
		\left\{
		P\in\PD(n)
		\;:\;
		P\text{ is an order-antichain distribution}
		\right\}.
		\]
	\end{definition}
	
	For \(s\in \{1,\eone,\etwo\}\) we have, \(U_n^s\in\mathcal A_n\), since
	all its coordinates are equal and the condition \(p_i\ne p_j\) never
	occurs. A distribution in \(\mathcal A_2\) distinct from the associated
	uniform distribution is	\(\left(0.7\eone+0.3\etwo,\;0.3\eone+0.7\etwo\right)\), whose two
	coordinates are incomparable. In contrast,
	\(\left(0.7\eone+0.7\etwo,\;0.3\eone+0.3\etwo\right)\notin\mathcal A_2\),
	since its second coordinate is strictly dominated by the first.
	\medskip
	
	A characterization of the fixed-point sets of hyperbolic
	negators on each mass stratum in terms of the class \(\mathcal A_n\)
	is established by our next theorem.
	
	\begin{theorem}
		\label{thm:fixed-set-realization}
		Let \(n\ge2\), \(s\in\{1,\eone,\etwo\}\), and let
		\(F\subseteq\PD(n)\) satisfy \(s(P)=s\) for every \(P\in F\).
		There exists a hyperbolic negator
		\(\Neg:\PD(n)\to\PD(n)\) such that
		\[
		\operatorname{Fix}_s(\Neg)=F
		\qquad \text{if and only if} \qquad
		U_n^s\in F\subseteq\mathcal A_n.
		\]
	\end{theorem}
	
	\begin{proof}
		Suppose first that
		\(F=\operatorname{Fix}_s(\Neg)\) for some hyperbolic negator
		\(\Neg\). By \eqref{eq:uniform-fixed-point}, \(U_n^s\in F\).
		Let \(P=(p_1,\ldots,p_n)\in F\), and suppose, to the contrary,
		that \(P\notin\mathcal A_n\). Then there exist \(i,j\) such that
		\(p_i\ne p_j\) and, without loss of generality,
		\(p_i\preceqD p_j\). By the negator property,
		\[
		\Neg(P)_j\preceqD\Neg(P)_i.
		\]
		Since \(\Neg(P)=P\), we obtain \(p_j\preceqD p_i\). By
		antisymmetry, \(p_i=p_j\), a contradiction. Hence
		\(F\subseteq\mathcal A_n\).
		
		Conversely, suppose that
		\(U_n^s\in F\subseteq\mathcal A_n\). Define
		\[
		\Neg_F(P)
		=
		\begin{cases}
			P, & P\in F,\\
			\Yager(P), & P\notin F.
		\end{cases}
		\]
		Both branches preserve hyperbolic mass. If \(P\notin F\), the
		order-reversing property follows from that of \(\Yager\). If
		\(P\in F\) and \(p_i\preceqD p_j\), then \(P\in\mathcal A_n\)
		implies \(p_i=p_j\), and therefore
		\[
		\Neg_F(P)_j=p_j=p_i=\Neg_F(P)_i \implies
		\Neg_F(P)_j\preceqD\Neg_F(P)_i.\] Thus \(\Neg_F\) is a
		hyperbolic negator.
		
		If \(P\in F\), then \(\Neg_F(P)=P\) by definition. Hence
		\(
		F\subseteq\operatorname{Fix}_s(\Neg_F)
		\).		
		For the reverse inclusion, let
		\(P\in\operatorname{Fix}_s(\Neg_F)\). If \(P\notin F\), then
		\(\Neg_F(P)=\Yager(P)\), and therefore \(\Yager(P)=P\).
		Thus
		\[
		p_i=\frac{s-p_i}{n-1}
		\implies
		p_i=\frac s n,\] for every \(i\). Hence
		\(P=U_n^s\in F\), a contradiction. Therefore \(P\in F\), and so
		\(\operatorname{Fix}_s(\Neg_F)\subseteq F\).
		Consequently,
		\[
		\operatorname{Fix}_s(\Neg_F)=F.
		\]
	\end{proof}
	
	\begin{remark}
		\label{rem:one-active-component-unique-fixed-point}
		In the strata \(s=\eone\) and \(s=\etwo\), any two coordinates are
		comparable. Hence, if \(P\in\mathcal A_n\) belongs to either of these
		strata, the order-antichain condition forces all its coordinates to be
		equal. The mass condition then gives
		\[
		\mathcal A_n\cap\{P\in\PD(n):s(P)=\eone\}
		=
		\{U_n^{\eone}\},
		\qquad
		\mathcal A_n\cap\{P\in\PD(n):s(P)=\etwo\}
		=
		\{U_n^{\etwo}\}.
		\]
		Therefore, by Theorem~\ref{thm:fixed-set-realization}, every
		hyperbolic negator satisfies
		\[
		\operatorname{Fix}_{\eone}(\Neg)=\{U_n^{\eone}\},
		\qquad
		\operatorname{Fix}_{\etwo}(\Neg)=\{U_n^{\etwo}\}.
		\]
		Consequently, a fixed point distinct from the associated uniform
		distribution may occur only in the stratum \(s=1\). In this stratum,
		if \(P=(p_1,\ldots,p_n)\in\mathcal A_n\) and \(p_i\ne p_j\), then
		either
		\[
		p_i^{(1)}<p_j^{(1)},\quad p_i^{(2)}>p_j^{(2)},
		\qquad\text{or}\qquad
		p_i^{(1)}>p_j^{(1)},\quad p_i^{(2)}<p_j^{(2)}.
		\]
	\end{remark}
	
	We now turn to conditions that recover uniqueness of the associated
	uniform distribution. Since the fixed-point structure in the strata
	\(s=\eone\) and \(s=\etwo\) is already determined, it remains to study
	the stratum \(s=1\). We first consider negators induced by admissible
	hyperbolic generators.
	
	For an invertible \(\lambda\in\Dpos\), define
	\[
	L_\lambda(\Phi)
	=
	\left\{
	z\in\Dint{0}{1}
	\;:\;
	\Phi(z)=\lambda z
	\right\}.
	\]
	
	\begin{theorem}
		\label{thm:level-set-fixed-point-criterion}
		Let \(n\ge2\), let
		\(\Phi:\Dint{0}{1}\to\Dpos\) be an admissible hyperbolic generator,
		and let \(\Neg_\Phi\) be the generated negator. If
		\(L_\lambda(\Phi)\) is totally ordered by \(\preceqD\) for every
		invertible \(\lambda\in\Dpos\), then \(\operatorname{Fix}_1(\Neg_\Phi)=\{U_n^1\}\).
	\end{theorem}
	
	\begin{proof}
		Let \(P=(p_1,\ldots,p_n)\in\operatorname{Fix}_1(\Neg_\Phi)\). By
		Theorem~\ref{thm:fixed-set-realization}, \(P\in\mathcal A_n\). Set
		\[
		\lambda_P
		=
		Z_1^\Phi(P)\eone+Z_2^\Phi(P)\etwo.
		\]
		Since \(s(P)=1\), both coefficients of \(\lambda_P\) are positive,
		and hence \(\lambda_P\) is invertible. For every
		\(k\in\{1,\ldots,n\}\), we have
		\[
		p_k
		=
		\Neg_\Phi(P)_k
		=
		\frac{\Phi(p_k)}{\lambda_P},
		\]
		and therefore \(\Phi(p_k)=\lambda_Pp_k\).
		Thus \(p_k\in L_{\lambda_P}(\Phi)\) for every \(k\). Since
		\(L_{\lambda_P}(\Phi)\) is totally ordered, any two coordinates of
		\(P\) are comparable. On the other hand, \(P\in\mathcal A_n\), so
		comparable coordinates must be equal. Hence
		\(p_1=\cdots=p_n\). Since \(s(P)=1\), it follows that
		\(P=U_n^1\). Therefore,
		\[
		\operatorname{Fix}_1(\Neg_\Phi)=\{U_n^1\}.
		\]
	\end{proof}
	
	\begin{remark}
		\label{rem:quotient-injectivity-criterion}
		A sufficient condition for the hypothesis of
		Theorem~\ref{thm:level-set-fixed-point-criterion} is that, for every
		invertible \(\lambda\in\Dpos\), the equation
		\(\Phi(z)=\lambda z\) has no noninvertible solution and
		\[
		z\longmapsto\frac{\Phi(z)}{z}
		\]
		is injective on the invertible elements of \(\Dint{0}{1}\).
		Indeed, under these assumptions, each \(L_\lambda(\Phi)\) contains
		at most one element, and is therefore totally ordered by
		\(\preceqD\).
	\end{remark}
	
\begin{figure}[H]
		\centering
		\includegraphics[width=0.8\textwidth]{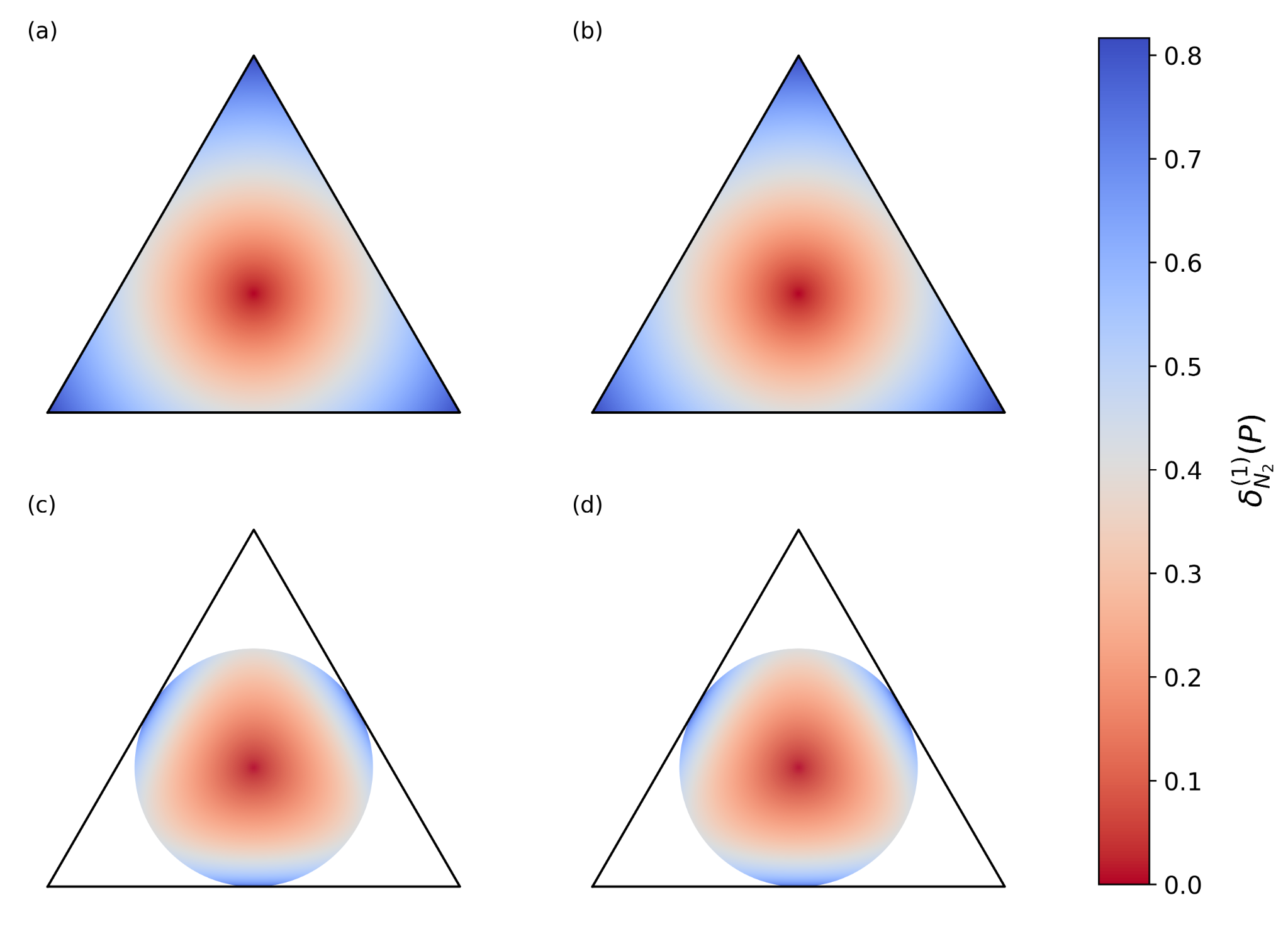}
		\caption{Visualization of the generated negator \(\Neg_2\) in
			Example~\ref{ex:quadratic-generated-fixed-point}. Panels
			\textup{(a)} and \textup{(b)} show the two idempotent components of
			the input distribution, and panels \textup{(c)} and \textup{(d)}
			those of its image under \(\Neg_2\). The color represents
			\(\delta_{\Neg_2}^{(1)}(P)\).}
		\label{fig:quadratic-generated-fixedset}
	\end{figure}

	\begin{example}
		\label{ex:quadratic-generated-fixed-point}
		Let \(n=3\) and consider the generator
		\[
		\Phi_2(z)=(1-z)^2,
		\qquad z\in\Dint{0}{1},
		\]
		with associated negator \(\Neg_2=\Neg_{\Phi_2}\). Since \(\Phi_2\) is
		order-reversing, \(\Neg_2\) is a hyperbolic negator.
		
		If \(z\) is noninvertible, the equation \(\Phi_2(z)=\lambda z\) is
		impossible, since a vanishing component of \(z\) corresponds to a
		component of \(\Phi_2(z)\) equal to \(1\). For invertible \(z\),
		\[
		\frac{\Phi_2(z)}{z}
		=
		g\bigl(z^{(1)}\bigr)\eone+g\bigl(z^{(2)}\bigr)\etwo,
		\qquad
		g(t)=\frac{(1-t)^2}{t}.
		\]
		Since \(g'(t)=1-t^{-2}<0\) on \((0,1)\), the function \(g\) is
		strictly decreasing on \((0,1]\), and hence the quotient is injective.
		Therefore, by Remark~\ref{rem:quotient-injectivity-criterion},
		\[
		\operatorname{Fix}_1(\Neg_2)=\{U_3^1\}.
		\]
		
		To visualize this fixed-point structure, we restrict the negator to the
		section
		\[
		\mathcal S=\{(x,Rx):x\in\Delta_3\},
		\qquad
		R(x_1,x_2,x_3)=(x_3,x_2,x_1).
		\]
		On this section, the two idempotent components of the Euclidean distance
		coincide. We use the scalar
		\[
		\delta_{\Neg_2}^{(1)}(P)
		=
		\inf_{Q\in\operatorname{Fix}_1(\Neg_2)}
		\bigl[d_E(P,Q)\bigr]^{(1)}
		=
		\bigl[d_E(P,U_3^1)\bigr]^{(1)}.
		\]
		In Figure~\ref{fig:quadratic-generated-fixedset}, each input
		distribution \(P\) is colored according to
		\(\delta_{\Neg_2}^{(1)}(P)\), and its image under \(\Neg_2\) retains
		the same color.
	\end{example}
	
\begin{figure}[h!]
		\centering
		\includegraphics[width=0.8\textwidth]{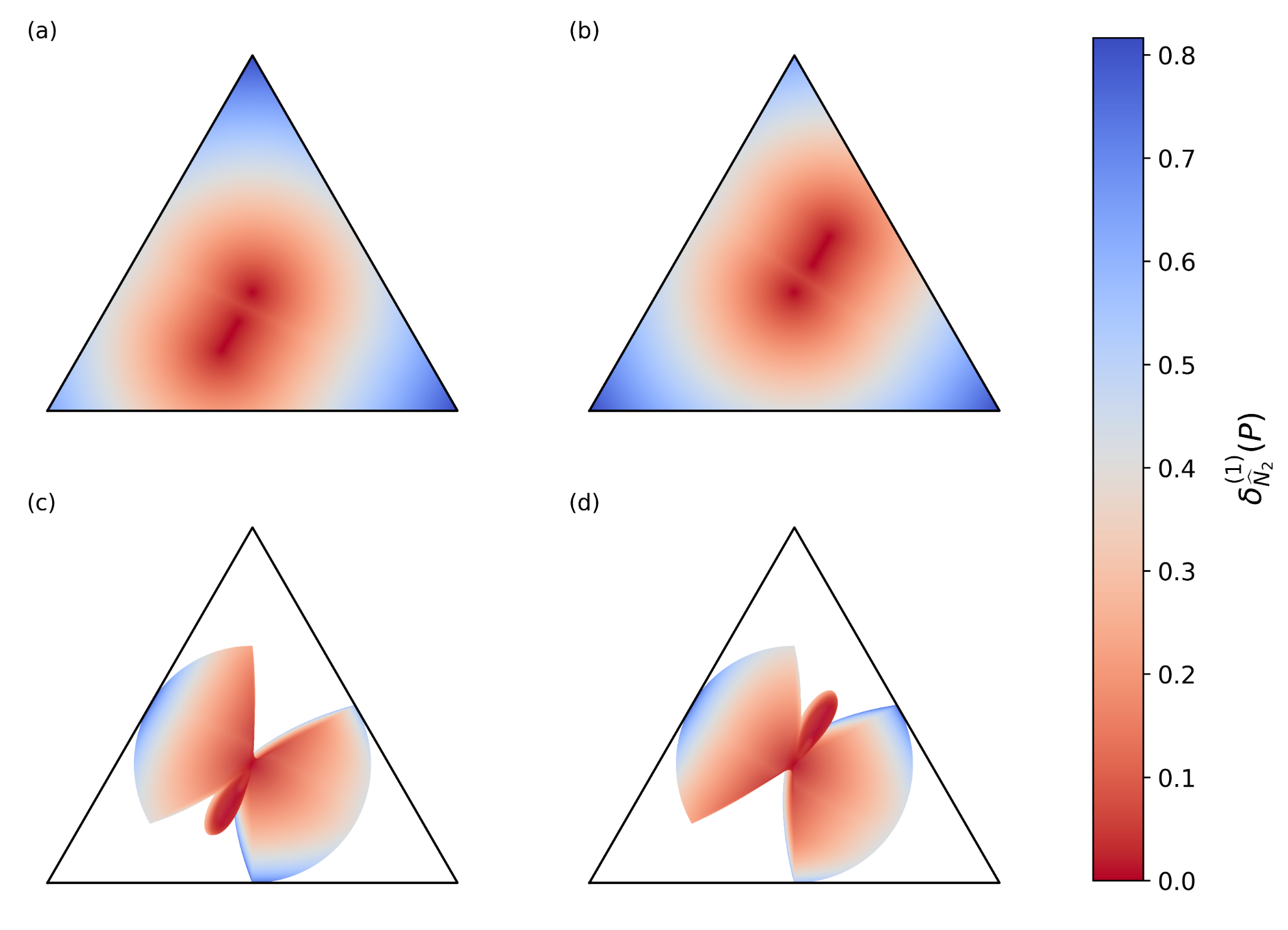}
		\caption{Visualization of the negator \(\widehat{\Neg}_2\) in
			Example~\ref{ex:continuous-negator-fixed-family}. Panels
			\textup{(a)} and \textup{(b)} show the two idempotent components of
			the input distribution, and panels \textup{(c)} and \textup{(d)}
			those of its image under \(\widehat{\Neg}_2\). The color represents
			\(\delta_{\widehat{\Neg}_2}^{(1)}(P)\).}
		\label{fig:continuous-negator-fixedset}
	\end{figure}

	\begin{example}
		\label{ex:continuous-negator-fixed-family}
		We now construct a hyperbolic negator with infinitely many
		additional fixed points.
		For \(t\in[0.08,0.16]\), let
		\(
		P_t=(0.33+t,\,0.34,\,0.33-t)\eone
		+(0.33-t,\,0.34,\,0.33+t)\etwo
		\)
		and set \(K=\{P_t\mid t\in[0.08,0.16]\}\). For every
		\(t\in[0.08,0.16]\), we have \(s(P_t)=1\), and the coordinates of
		\(P_t\) are pairwise distinct and incomparable. Hence \(P_t\) is an
		order-antichain distribution.
		
		For \(r\in\mathbb R\), write \(r_+=\max\{r,0\}\). Given
		\(P=(p_1,p_2,p_3)\in\PD(3)\), define
		\[
		\alpha(P)=
		\min_{1\le i<j\le3}
		\left[-(p_i^{(1)}-p_j^{(1)})(p_i^{(2)}-p_j^{(2)})\right]_+ .
		\]
		Thus \(\alpha(P)>0\) precisely when the three coordinates of \(P\)
		are pairwise distinct and pairwise incomparable. Set
		\(
		\tau(P)=\frac{p_1^{(1)}-p_1^{(2)}}2
		\)
		and
		\[
		\rho(P)=
		(p_2^{(1)}-0.34)^2+(p_2^{(2)}-0.34)^2
		+(p_1^{(1)}+p_1^{(2)}-0.66)^2
		+(0.08-\tau(P))_+^2+(\tau(P)-0.16)_+^2 .
		\]
		Then \(\rho(P)=0\) if and only if \(P\in K\).
		Define
		\[
		w(P)=\frac{\rho(P)}{\rho(P)+\alpha(P)},
		\qquad
		\widehat{\Neg}_2(P)=(1-w(P))P+w(P)\Neg_2(P).
		\]
		The denominator is positive because \(\rho(P)=0\) implies
		\(P\in K\) and hence \(\alpha(P)>0\). Therefore
		\(\widehat{\Neg}_2\) is continuous, preserves hyperbolic mass, and maps
		\(\PD(3)\) into itself.
		
		If \(p_i\preceqD p_j\) for some \(i\ne j\), then \(\alpha(P)=0\), so
		\(w(P)=1\) and \(\widehat{\Neg}_2(P)=\Neg_2(P)\). Hence
		\(\widehat{\Neg}_2\) is a hyperbolic negator.
		
		Moreover, \(w(P_t)=0\) for every \(P_t\in K\), and therefore
		\(\widehat{\Neg}_2(P_t)=P_t\). Conversely, if
		\(\widehat{\Neg}_2(P)=P\) and \(P\notin K\), then \(w(P)>0\), so
		\[
		P=(1-w(P))P+w(P)\Neg_2(P)
		\quad\to\quad
		\Neg_2(P)=P.
		\]
		Example~\ref{ex:quadratic-generated-fixed-point} gives \(P=U_3^1\).
		Consequently,
		\[
		\operatorname{Fix}_1(\widehat{\Neg}_2)=K\cup\{U_3^1\}.
		\]
		
		Figure~\ref{fig:continuous-negator-fixedset} illustrates the action of
		\(\widehat{\Neg}_2\) on the same section \(\mathcal S\). We use the
		scalar
		\[
		\delta_{\widehat{\Neg}_2}^{(1)}(P)
		=
		\inf_{Q\in\operatorname{Fix}_1(\widehat{\Neg}_2)}
		\bigl[d_E(P,Q)\bigr]^{(1)}.
		\]
		Each input distribution \(P\) is colored according to
		\(\delta_{\widehat{\Neg}_2}^{(1)}(P)\), and its image under
		\(\widehat{\Neg}_2\) retains the same color.
	\end{example}
	
	Quasi-contractive mappings have been extensively studied in cone-valued
	metric spaces, including nonnormal cones and nonlinear conditions governed
	by comparison mappings on the cone
	\citep{ArandelovicKeckic2011,GajicRakocevic2012}. Related conditions arise
	in the theory of quasi-nonexpansive mappings in cone metric spaces
	\citep{Saluja2013}. In the present setting, these ideas can be expressed
	directly in terms of hyperbolic-valued disequilibrium. Since the associated
	uniform distribution is fixed by every hyperbolic negator, a comparison
	map lying strictly below the identity away from zero provides a sufficient
	condition for its uniqueness.

	\FloatBarrier

	\begin{theorem}
		\label{thm:strict-quasi-nonexpansive-uniqueness}
		Let \(d\) be a hyperbolic-valued metric on each mass stratum of
		\(\PD(n)\), with associated disequilibrium \(D_d\), and let
		\(\Neg:\PD(n)\to\PD(n)\) be a hyperbolic negator. Suppose that there
		exists a map \(\varphi:\Dpos\to\Dpos\) such that, for each $z\in\Dpos$, $\varphi(z)\preceqD z$, with equality if and only if \(z=0\), and
		\[
		D_d\bigl(\Neg(P)\bigr)
		\preceqD
		\varphi\bigl(D_d(P)\bigr)
		\]
		for every \(P\in\PD(n)\). Then, for $s\in\{1,\eone,\etwo\}$, we have
		\[
		\operatorname{Fix}_s(\Neg)=\{U_n^s\}.
		\]
	\end{theorem}
	
	\begin{proof}
		Let \(P\in\operatorname{Fix}_s(\Neg)\). Since \(\Neg(P)=P\),
		\[
		D_d(P)
		=
		D_d\bigl(\Neg(P)\bigr)
		\preceqD
		\varphi\bigl(D_d(P)\bigr)
		\preceqD
		D_d(P).
		\]
		By antisymmetry,
		\[
		\varphi\bigl(D_d(P)\bigr)=D_d(P),
		\]
		and therefore \(D_d(P)=0\). Since \(D_d\) is induced by a
		hyperbolic-valued metric, this implies
		\(P=U_n^{s(P)}\). Together with
		\eqref{eq:uniform-fixed-point}, we obtain
		\[
		\operatorname{Fix}_s(\Neg)=\{U_n^s\}.
		\]
	\end{proof}
	
	The preceding result concerns uniqueness of the	fixed point. In the next section, we strengthen the comparison conditions
	to study the iterates of a hyperbolic negator and establish conditions
	under which $D_d\bigl(\Neg^m(P)\bigr)\to 0$.

	\section{Decay of hyperbolic-valued disequilibrium}
	\label{sec:disequilibrium-decay}

	\citet{BoryReyes2026} developed a Yager-based majorization framework
	for a broad class of generated hyperbolic negators, establishing
	entropy increase and quantitative convergence toward the associated
	uniform distribution. In the classical setting,
	\citet{BrunoVaccaro2026} obtained related quantitative bounds for
	entropy growth under iterated independent negation. Here, we describe this behavior
	in terms of hyperbolic-valued disequilibrium, derive
	decay estimates for \(D_E\), \(D_h\), and
	\(D_{\operatorname{JS}}\), and develop an independent
	Lipschitz--spectral approach for generated negators.

	The preceding section provided a comparison condition ensuring uniqueness
	of the associated uniform distribution as a fixed point. To study the
	iterates of a hyperbolic negator, this condition is strengthen by requiring
	the comparison map itself to converge to zero under iteration. The
	following criterion applies to any hyperbolic-valued disequilibrium.
	
	\begin{theorem}
		\label{thm:disequilibrium-decay-criterion}
		Let \(d\) be a hyperbolic-valued metric on each mass stratum of
		\(\PD(n)\), with associated disequilibrium \(D_d\), and let
		\(\Neg:\PD(n)\to\PD(n)\) be a hyperbolic negator. Suppose that there
		exists an order-preserving map
		\(\varphi:\Dpos\to\Dpos\) such that $\varphi^m(z)\to 0$	for every \(z\in\Dpos\), and
		\[
		D_d\bigl(\Neg(Q)\bigr)
		\preceqD
		\varphi\bigl(D_d(Q)\bigr)
		\]
		for every \(Q\in\PD(n)\). Then, for every \(P\in\PD(n)\), $$D_d\bigl(\Neg^m(P)\bigr)\to 0.$$
	\end{theorem}
	
	\begin{proof}
		Since \(\Neg\) preserves hyperbolic mass, $s\bigl(\Neg^m(P)\bigr)=s(P)$.
		Applying the assumed inequality successively and using the
		order-preserving property of \(\varphi\), we obtain
		\[
		D_d\bigl(\Neg^m(P)\bigr)
		\preceqD
		\varphi^m\bigl(D_d(P)\bigr).
		\]
		Since the right-hand side converges to \(0\), so does
		\(D_d(\Neg^m(P))\).
	\end{proof}
	
	\begin{remark}
		\label{rem:linear-comparison-map}
		If \(\varphi:\Dpos\to\Dpos\) is a positive real-linear map, then,
		with respect to the idempotent basis, it is represented by a
		non-negative matrix
		\[
		A=
		\begin{pmatrix}
			a_{11}&a_{12}\\
			a_{21}&a_{22}
		\end{pmatrix}.
		\]
		In this case, $\varphi^m(z)\to 0$ for every $z\in\Dpos$ 
		if and only if \(\rho(A)<1\). Thus the spectral-radius condition
		provides the linear form of
		Theorem~\ref{thm:disequilibrium-decay-criterion}.
	\end{remark}
	
	Next, two mechanisms leading the previous comparative estimate are outlined: One based on domination by Yager's negator and the other on Lipschitz
	bounds for generated negators.
	
	\subsection{Yager-dominated negators}
	\label{subsec:yager-dominated-negators}
	
	The majorization comparison in
	Proposition~\ref{prop:yager-intermediate-majorization}, together with the
	hyperbolic Schur-convexity established in
	Proposition~\ref{prop:uniform-metrics-schur-convex}, yields direct decay
	estimates for \(D_E\), \(D_h\), and \(D_{\operatorname{JS}}\).
	We first quantify the action of Yager's negator on these three
	disequilibria.
	
	Define
	\[
	\theta_n^{\operatorname{JS}}
	=
	\inf\left\{
	\theta\ge0:
	D_{\operatorname{JS}}\bigl(\Yager(P)\bigr)
	\preceqD
	\theta\,D_{\operatorname{JS}}(P)
	\text{ for every }P\in\PD(n)
	\right\}.
	\]
	By Propositions~\ref{prop:yager-intermediate-majorization} and
	\ref{prop:uniform-metrics-schur-convex},
	$
	D_{\operatorname{JS}}\bigl(\Yager(P)\bigr)
	\preceqD
	D_{\operatorname{JS}}(P)
	$,
	so the set in the preceding definition is nonempty.
	
	\begin{proposition}
		\label{prop:yager-disequilibrium-decay}
		Let \(n>2\). Then
		\(\theta_n^{\operatorname{JS}}\in[0,1)\), and, for every
		\(P\in\PD(n)\),
		\[
		D_h\bigl(\Yager(P)\bigr)
		=
		\frac{1}{n-1}D_h(P),
		\]
		\[
		D_E\bigl(\Yager(P)\bigr)
		=
		\frac{1}{n-1}D_E(P),
		\]
		and
		\[
		D_{\operatorname{JS}}\bigl(\Yager(P)\bigr)
		\preceqD
		\theta_n^{\operatorname{JS}}\,
		D_{\operatorname{JS}}(P).
		\]
	\end{proposition}
	
	\begin{proof}
		For every \(j\in\{1,\ldots,n\}\),
		\[
		\Yager(P)_j-\frac{s(P)}{n}
		=
		-\frac{1}{n-1}
		\left(
		p_j-\frac{s(P)}{n}
		\right).
		\]
		The identities for \(D_h\) and \(D_E\) follow directly from their
		definitions.
		\medskip
		
		Finally, let us consider \(D_{\operatorname{JS}}\). Fix
		\(\ell\in\{1,2\}\) such that \(s^{(\ell)}(P)=1\). Application of
		Propositions~\ref{prop:yager-intermediate-majorization} and \ref{prop:uniform-metrics-schur-convex}, give
		\[
		\left[
		D_{\operatorname{JS}}\bigl(\Yager(P)\bigr)
		\right]^{(\ell)}
		<
		\left[D_{\operatorname{JS}}(P)\right]^{(\ell)}.
		\]
		
		It remains to obtain a factor independent of \(P\). For
		$
		P^{(\ell)}
		\ne
		\left(\frac1n,\ldots,\frac1n\right)
		$,
		define
		\[
		\eta_n^{\operatorname{JS}}\bigl(P^{(\ell)}\bigr)
		=
		\frac{
			\left[
			D_{\operatorname{JS}}\bigl(\Yager(P)\bigr)
			\right]^{(\ell)}
		}{
			\left[D_{\operatorname{JS}}(P)\right]^{(\ell)}
		}.
		\]
		This quantity depends only on \(P^{(\ell)}\), and hence defines a
		function on
		$
		\Delta_n\setminus
		\{(\frac1n,\ldots,\frac1n)\}
		$.
		
		Let \(\psi_n\) be the function used in the proof of
		Proposition~\ref{prop:uniform-metrics-schur-convex}. Then
		$
		\psi_n\left(\frac1n\right)
		=
		\psi_n'\left(\frac1n\right)
		=
		0
		$,
		and
		$
		\psi_n''\left(\frac1n\right)
		=
		\frac n4
		$.
		Set
		$
		h_j=p_j^{(\ell)}-\frac1n
		$.
		For \(P^{(\ell)}\) sufficiently close to
		\((\frac1n,\ldots,\frac1n)\), the Taylor theorem ensures that for each \(j\),
		\[
		\psi_n\left(p_j^{(\ell)}\right)
		=
		\psi_n\left(\frac1n+h_j\right)
		=
		\psi_n\left(\frac1n\right)
		+
		\psi_n'\left(\frac1n\right)h_j
		+
		\frac12\psi_n''(\xi_j)h_j^2
		=
		\frac12\psi_n''(\xi_j)h_j^2,
		\]
		and some \(\xi_j\) between \(\frac1n\) and
		\(\frac1n+h_j\). Hence, whenever
		$
		P^{(\ell)}
		\ne
		\left(\frac1n,\ldots,\frac1n\right)
		$, we have
		\[
		\frac{
			\left([D_{\operatorname{JS}}(P)]^{(\ell)}\right)^2
		}{
			\sum_{j=1}^n h_j^2
		}
		=
		\frac12
		\sum_{j=1}^n
		\psi_n''(\xi_j)
		\frac{h_j^2}{\sum_{k=1}^n h_k^2}.
		\]
		Since
		$
		\psi_n''(t)=\frac{1}{2n\,t(t+1/n)}
		$
		is continuous at \(\frac1n\),
		$
		\psi_n''\left(\frac1n\right)=\frac n4
		$,
		and
		$
		|\xi_j-\frac1n|\le |h_j|
		$,
		it follows that
		\begin{equation}\label{eq:Tay_Yager}
			\lim_{
				P^{(\ell)}
				\to
				(\frac1n,\ldots,\frac1n)
			}
			\frac{
				\left([D_{\operatorname{JS}}(P)]^{(\ell)}\right)^2
			}{
				\sum_{j=1}^n
				\left(p_j^{(\ell)}-\frac1n\right)^2
			}
			=
			\frac n8.
		\end{equation}
		
		Moreover,
		\[
		\sum_{j=1}^n
		\left(
		\bigl(\Yager(P)\bigr)_j^{(\ell)}-\frac1n
		\right)^2
		=
		\frac1{(n-1)^2}
		\sum_{j=1}^n
		\left(p_j^{(\ell)}-\frac1n\right)^2.
		\]
		Since
		$
		P^{(\ell)}
		\to
		(\frac1n,\ldots,\frac1n)
		$
		implies
		$
		(\Yager(P))^{(\ell)}
		\to
		(\frac1n,\ldots,\frac1n)
		$,
		Equation~\eqref{eq:Tay_Yager}, applied to \(P\) and
		\(\Yager(P)\), yields
		\[
		\lim_{
			P^{(\ell)}
			\to
			(\frac1n,\ldots,\frac1n)
		}
		\left(
		\eta_n^{\operatorname{JS}}\bigl(P^{(\ell)}\bigr)
		\right)^2
		=
		\frac n8
		\frac1{(n-1)^2}
		\frac8n
		=
		\frac1{(n-1)^2}.
		\]
		Therefore,
		\[
		\lim_{
			P^{(\ell)}
			\to
			(\frac1n,\ldots,\frac1n)
		}
		\eta_n^{\operatorname{JS}}\bigl(P^{(\ell)}\bigr)
		=
		\frac1{n-1}.
		\]
		
		Defining
		$
		\eta_n^{\operatorname{JS}}
		\left(\frac1n,\ldots,\frac1n\right)
		=
		\frac1{n-1}
		$
		gives a continuous extension of
		\(\eta_n^{\operatorname{JS}}\) to \(\Delta_n\). The strict
		Schur-convexity argument above gives
		$
		\eta_n^{\operatorname{JS}}(P^{(\ell)})<1
		$
		whenever
		$
		P^{(\ell)}
		\ne
		(\frac1n,\ldots,\frac1n)
		$,
		while \(n>2\) gives
		$
		\eta_n^{\operatorname{JS}}
		(\frac1n,\ldots,\frac1n)
		=
		1/(n-1)<1
		$.
		Hence, by compactness,
		\[
		M_n
		=
		\max_{x\in\Delta_n}
		\eta_n^{\operatorname{JS}}(x)
		<1.
		\]
		Since \(M_n\) is admissible in the definition of
		\(\theta_n^{\operatorname{JS}}\), it follows that
		\[
		0\le
		\theta_n^{\operatorname{JS}}
		\le
		M_n
		<1.
		\]
	\end{proof}
	
	The restriction \(n>2\) is essential. When \(n=2\), Yager's negator
	exchanges the two coordinates in each active coefficient. Since
	\(D_E\), \(D_h\), and \(D_{\operatorname{JS}}\) are invariant under
	permutations,
	\[
	D_\nu\bigl(\Yager(P)\bigr)=D_\nu(P),
	\qquad
	\nu\in\{E,h,\operatorname{JS}\}.
	\]	
	
	\begin{theorem}
		\label{thm:yager-dominated-disequilibrium-decay}
		Let \(n>2\), and let \(\Neg:\PD(n)\to\PD(n)\) be a hyperbolic negator
		satisfying
		$
		\Yager(P)\Dmaj\Neg(P)
		$,
		for every \(P\in\PD(n)\). Then, for every \(P\in\PD(n)\),
		\[
		D_h\bigl(\Neg^m(P)\bigr)\to0,
		\qquad
		D_E\bigl(\Neg^m(P)\bigr)\to0,
		\qquad \text{and} \qquad
		D_{\operatorname{JS}}\bigl(\Neg^m(P)\bigr)\to0.
		\]
	\end{theorem}
	
	\begin{proof}
		By Proposition~\ref{prop:uniform-metrics-schur-convex}, the relation
		$
		\Yager(Q)\Dmaj\Neg(Q)
		$
		implies
		\[
		D_\nu\bigl(\Neg(Q)\bigr)
		\preceqD
		D_\nu\bigl(\Yager(Q)\bigr),
		\qquad
		\nu\in\{h,E,\operatorname{JS}\}.
		\]
		Therefore 
		\[
		D_h\bigl(\Neg(Q)\bigr)
		\preceqD
		\frac{1}{n-1}D_h(Q),
		\qquad
		D_E\bigl(\Neg(Q)\bigr)
		\preceqD
		\frac{1}{n-1}D_E(Q),
		\]
		by Proposition~\ref{prop:yager-disequilibrium-decay} and thus
		\[
		D_{\operatorname{JS}}\bigl(\Neg(Q)\bigr)
		\preceqD
		\theta_n^{\operatorname{JS}}
		D_{\operatorname{JS}}(Q).
		\]
		Since
		$
		1/(n-1)<1
		$
		and
		$
		\theta_n^{\operatorname{JS}}<1
		$,
		Theorem~\ref{thm:disequilibrium-decay-criterion}, applied with
		$
		\varphi(z)=z/(n-1)
		$
		for \(D_h\) and \(D_E\), and with
		$
		\varphi(z)=\theta_n^{\operatorname{JS}}z
		$
		for \(D_{\operatorname{JS}}\), yields the result.
	\end{proof}
	
\begin{figure}[htbp]
		\centering
		\includegraphics[width=\textwidth]
		{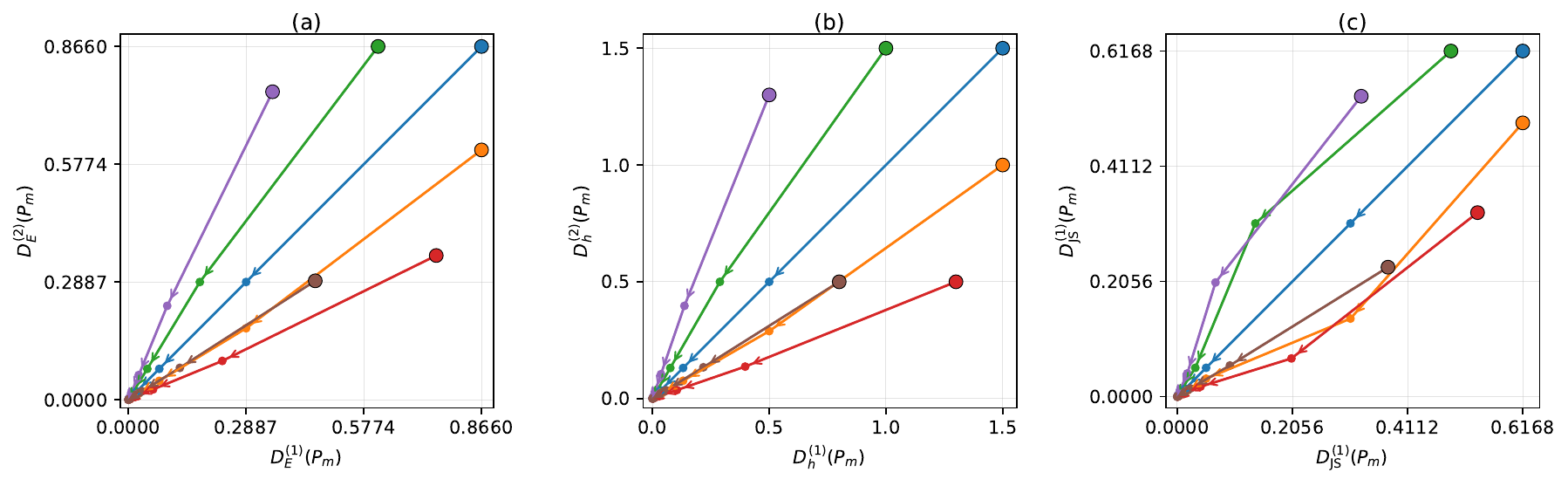}
		\caption{Decay of the disequilibrium associated with the iterates of
			\(\Neg_\alpha\) from
			Example~\ref{ex:yager-dominated-power-generated-negator}, with
			\(n=4\) and \(\alpha=0.8\). Panels \textup{(a)}, \textup{(b)}, and
			\textup{(c)} display \(D_E\), \(D_h\), and
			\(D_{\operatorname{JS}}\), respectively.}
		\label{fig:yager-dominated-negator-metric-orbits}
	\end{figure}

	\begin{example}
		\label{ex:yager-dominated-power-generated-negator}
		Let \(n>2\), \(0<\alpha\le1\), and define
		\[
		\Phi_\alpha(z)=(1-z)^\alpha,
		\qquad z\in\Dint{0}{1}.
		\]
		This is an admissible hyperbolic generator, and we denote the
		corresponding negator by
		$
		\Neg_\alpha=\Neg_{\Phi_\alpha}
		$.
		For \(\alpha=1\), one has \(\Neg_1=\Yager\).
		
		For every active component,
		\[
		\left[\Phi_\alpha(p_j)\right]^{(\ell)}
		=
		\left(1-p_j^{(\ell)}\right)^\alpha,
		\]
		and, whenever \(p_j^{(\ell)}<1\),
		\[
		\frac{
			\left[\Phi_\alpha(p_j)\right]^{(\ell)}
		}{
			1-p_j^{(\ell)}
		}
		=
		\left(1-p_j^{(\ell)}\right)^{\alpha-1}.
		\]
		Since \(0<\alpha\le1\), the map
		\(t\mapsto(1-t)^\alpha\) is non-increasing on \([0,1]\), whereas
		\(t\mapsto(1-t)^{\alpha-1}\) is non-decreasing on \([0,1)\).
		Moreover,
		$
		p_j^{(\ell)}=1
		$
		implies
		$
		[\Phi_\alpha(p_j)]^{(\ell)}=0.
		$
		Thus the hypotheses of
		Theorem~\ref{thm:generated-negator-yager-dominated} hold, and
		\[
		\Yager(P)\Dmaj\Neg_\alpha(P),
		\qquad P\in\PD(n).
		\]
		Therefore, by
		Theorem~\ref{thm:yager-dominated-disequilibrium-decay}, for every
		\(P\in\PD(n)\),
		\[
		D_h\bigl(\Neg_\alpha^m(P)\bigr)\longrightarrow0,
		\qquad
		D_E\bigl(\Neg_\alpha^m(P)\bigr)\longrightarrow0,
		\qquad
		D_{\operatorname{JS}}\bigl(\Neg_\alpha^m(P)\bigr)\longrightarrow0.
		\]
		Figure~\ref{fig:yager-dominated-negator-metric-orbits} illustrates this
		decay for \(n=4\) and \(\alpha=0.8\).
	\end{example}
	
	\subsection{Lipschitz conditions for generated negators}
	\label{subsec:lipschitz-generated-negators}
	
	We now describe the second mechanism for the decay of hyperbolic-valued
	disequilibrium, based on Lipschitz-type estimates for the admissible
	generator defining the negator. Related notions of
	\(\mathbb D\)-Lipschitz continuity and their use in the convergence
	analysis of iterative schemes in hyperbolic Hilbert spaces have been
	studied by \citet{AliKumar2026}.
	
	Let \(\Phi:\Dint{0}{1}\to\Dpos\) be an admissible hyperbolic generator.
	Assume that there exist constants \(L_{\ell r}\ge0\),
	\(\ell,r\in\{1,2\}\), such that
	\begin{equation}
		\label{eq:generator-lipschitz-condition}
		\left|
		\Phi^{(\ell)}(z)-\Phi^{(\ell)}(w)
		\right|
		\le
		L_{\ell1}|z^{(1)}-w^{(1)}|
		+
		L_{\ell2}|z^{(2)}-w^{(2)}|
	\end{equation}
	for every \(z,w\in\Dint{0}{1}\) and every \(\ell\in\{1,2\}\).
	
	Fix \(s\in\{1,\eone,\etwo\}\). For every active component \(\ell\),
	assume that there exists a constant \(M_{\ell,s}>0\) such that
	\begin{equation}
		\label{eq:normalizing-factor-lower-bound}
		Z_\ell^\Phi(Q)\ge M_{\ell,s}
	\end{equation}
	for every \(Q\in\PD(n)\) with \(s(Q)=s\).
	
	We shall use \(\psi_n\) introduced in the proof of
	Proposition~\ref{prop:uniform-metrics-schur-convex} again. For every
	\(P\in\PD(n)\) and every active component \(\ell\),
	\[
	\left(
	\left[D_{\operatorname{JS}}(P)\right]^{(\ell)}
	\right)^2
	=
	\sum_{j=1}^{n}\psi_n\left(p_j^{(\ell)}\right).
	\]
	The Taylor calculation leading to Equation~\eqref{eq:Tay_Yager} gives
	\[
	\lim_{t\to1/n}
	\frac{\psi_n(t)}{(t-1/n)^2}
	=
	\frac n8.
	\]
	Accordingly, define
	\[
	\chi_n(t)
	=
	\begin{cases}
		\displaystyle
		\frac{\psi_n(t)}{(t-1/n)^2},
		& t\ne\frac1n,\\[10pt]
		\displaystyle
		\frac n8,
		& t=\frac1n.
	\end{cases}
	\]
	Thus \(\chi_n\) is positive and continuous on \([0,1]\). Set \(\mu_n = \min_{0\le t\le1}\chi_n(t)\), \(\nu_n=\max_{0\le t\le1}\chi_n(t)\).
	Hence
	$
	0<\mu_n\le\nu_n<\infty
	$.
	Define
	$
	\eta_n
	=
	\sqrt{\frac{\nu_n}{\mu_n}}
	$.
	
	Now define the matrices
	\(A_{h,s}=\left(a_{\ell r}^{h,s}\right)\),
	\(A_{E,s}=\left(a_{\ell r}^{E,s}\right)\), and
	\(A_{\operatorname{JS},s}
	=\left(a_{\ell r}^{\operatorname{JS},s}\right)\) by
	\[
	a_{\ell r}^{h,s}
	=
	\begin{cases}
		\displaystyle
		\frac{2L_{\ell r}}{M_{\ell,s}},
		& \text{if }s^{(\ell)}=1,\\[8pt]
		0,
		& \text{if }s^{(\ell)}=0,
	\end{cases}
	\quad
	a_{\ell r}^{E,s}
	=
	\begin{cases}
		\displaystyle
		\frac{L_{\ell r}}{M_{\ell,s}},
		& \text{if }s^{(\ell)}=1,\\[8pt]
		0,
		& \text{if }s^{(\ell)}=0,
	\end{cases}
	\quad
	\text{and}
	\quad
	a_{\ell r}^{\operatorname{JS},s}
	=
	\begin{cases}
		\displaystyle
		\eta_n\,\frac{L_{\ell r}}{M_{\ell,s}},
		& \text{if }s^{(\ell)}=1,\\[8pt]
		0,
		& \text{if }s^{(\ell)}=0.
	\end{cases}
	\]
	
	For a non-negative matrix \(A\), let
	\(\mathcal A_A:\Dpos\to\Dpos\) denote the positive real-linear map
	represented by \(A\) with respect to the idempotent basis.
	
	\begin{proposition} 
		\label{prop:lipschitz-generated-negators-disequilibrium-decay}
		Let \(n\ge2\), let \(\Phi:\Dint{0}{1}\to\Dpos\) be an admissible
		hyperbolic generator satisfying
		\eqref{eq:generator-lipschitz-condition}, and let
		\(\Neg_\Phi\) be the generated negator. Fix
		\(s\in\{1,\eone,\etwo\}\), and suppose that
		\eqref{eq:normalizing-factor-lower-bound} holds for every active
		component \(\ell\). Then:
		
		\begin{enumerate}[label=\textup{(\roman*)}]
			\item If \(\rho(A_{h,s})<1\), then, for every \(P\in\PD(n)\)
			with \(s(P)=s\),
			\[
			D_h\bigl(\Neg_\Phi^m(P)\bigr)\longrightarrow0.
			\]
			
			\item If \(\rho(A_{E,s})<1\), then, for every \(P\in\PD(n)\)
			with \(s(P)=s\),
			\[
			D_E\bigl(\Neg_\Phi^m(P)\bigr)\longrightarrow0.
			\]
			
			\item If \(\rho(A_{\operatorname{JS},s})<1\), then, for every
			\(P\in\PD(n)\) with \(s(P)=s\),
			\[
			D_{\operatorname{JS}}
			\bigl(\Neg_\Phi^m(P)\bigr)\longrightarrow0.
			\]
		\end{enumerate}
	\end{proposition}
	
	\begin{proof}
		Let \(Q=(q_1,\ldots,q_n)\in\PD(n)\) with \(s(Q)=s\), and write
		\[
		R=\Neg_\Phi(Q)=(r_1,\ldots,r_n).
		\]
		Fix an active component \(\ell\). Then
		\[
		r_j^{(\ell)}
		=
		\frac{\Phi^{(\ell)}(q_j)}{Z_\ell^\Phi(Q)},
		\qquad j=1,\ldots,n.
		\]
		Since \(s^{(\ell)}=1\),
		\[
		\frac1n
		=
		\frac{1}{Z_\ell^\Phi(Q)}
		\frac1n
		\sum_{k=1}^{n}\Phi^{(\ell)}(q_k).
		\]
		Hence
		\[
		r_j^{(\ell)}-\frac1n
		=
		\frac{
			\Phi^{(\ell)}(q_j)
			-
			\frac1n\sum_{k=1}^{n}\Phi^{(\ell)}(q_k)
		}{
			Z_\ell^\Phi(Q)
		}.
		\]
		
		We first prove the estimate for \(D_h\). For every real number \(c\),
		\[
		\sum_{j=1}^{n}
		\left|
		\Phi^{(\ell)}(q_j)
		-
		\frac1n\sum_{k=1}^{n}\Phi^{(\ell)}(q_k)
		\right|
		\le
		2\sum_{j=1}^{n}
		\left|
		\Phi^{(\ell)}(q_j)-c
		\right|.
		\]
		Taking \(c=\Phi^{(\ell)}(s/n)\), and using
		\eqref{eq:generator-lipschitz-condition} together with
		\eqref{eq:normalizing-factor-lower-bound}, we obtain
		\[
		\begin{aligned}
			\sum_{j=1}^{n}
			\left|
			r_j^{(\ell)}-\frac1n
			\right|
			&\le
			\frac{2}{M_{\ell,s}}
			\sum_{j=1}^{n}
			\left|
			\Phi^{(\ell)}(q_j)-\Phi^{(\ell)}(s/n)
			\right|\\
			&\le
			\sum_{r=1}^{2}
			\frac{2L_{\ell r}}{M_{\ell,s}}
			\sum_{j=1}^{n}
			\left|
			q_j^{(r)}-\frac{s^{(r)}}{n}
			\right|.
		\end{aligned}
		\]
		If \(s^{(\ell)}=0\), then
		$
		q_j^{(\ell)}=r_j^{(\ell)}=0
		$
		for every \(j\), and the corresponding row of \(A_{h,s}\) is zero.
		Therefore
		\begin{equation}
			\label{eq:lipschitz-Dh-estimate}
			D_h(R)
			\preceqD
			\mathcal A_{A_{h,s}}
			\bigl(D_h(Q)\bigr).
		\end{equation}
		
		We now prove the estimate for \(D_E\). Since the arithmetic mean
		minimizes the quadratic deviation,
		\[
		\left(
		\sum_{j=1}^{n}
		\left|
		\Phi^{(\ell)}(q_j)
		-
		\frac1n\sum_{k=1}^{n}\Phi^{(\ell)}(q_k)
		\right|^2
		\right)^{1/2}
		\le
		\left(
		\sum_{j=1}^{n}
		\left|
		\Phi^{(\ell)}(q_j)-\Phi^{(\ell)}(s/n)
		\right|^2
		\right)^{1/2}.
		\]
		Using again
		\eqref{eq:generator-lipschitz-condition} and
		\eqref{eq:normalizing-factor-lower-bound}, we obtain
		\[
		\begin{aligned}
			\left(
			\sum_{j=1}^{n}
			\left|
			r_j^{(\ell)}-\frac1n
			\right|^2
			\right)^{1/2}
			&\le
			\frac1{M_{\ell,s}}
			\left(
			\sum_{j=1}^{n}
			\left|
			\Phi^{(\ell)}(q_j)-\Phi^{(\ell)}(s/n)
			\right|^2
			\right)^{1/2}\\
			&\le
			\sum_{r=1}^{2}
			\frac{L_{\ell r}}{M_{\ell,s}}
			\left(
			\sum_{j=1}^{n}
			\left|
			q_j^{(r)}-\frac{s^{(r)}}{n}
			\right|^2
			\right)^{1/2}.
		\end{aligned}
		\]
		Inactive components vanish as above. Therefore
		\begin{equation}
			\label{eq:lipschitz-DE-estimate}
			D_E(R)
			\preceqD
			\mathcal A_{A_{E,s}}
			\bigl(D_E(Q)\bigr).
		\end{equation}
		
		It remains to prove the estimate for \(D_{\operatorname{JS}}\).
		From the definitions of \(\mu_n\) and \(\nu_n\), for every active
		component \(\ell\),
		\[
		\begin{aligned}
			\left[D_{\operatorname{JS}}(R)\right]^{(\ell)}
			&\le
			\sqrt{\nu_n}
			\left(
			\sum_{j=1}^{n}
			\left|
			r_j^{(\ell)}-\frac1n
			\right|^2
			\right)^{1/2}\\
			&\le
			\sum_{r=1}^{2}
			\sqrt{\nu_n}\,
			\frac{L_{\ell r}}{M_{\ell,s}}
			\left(
			\sum_{j=1}^{n}
			\left|
			q_j^{(r)}-\frac{s^{(r)}}{n}
			\right|^2
			\right)^{1/2}\\
			&\le
			\sum_{r=1}^{2}
			\eta_n\,
			\frac{L_{\ell r}}{M_{\ell,s}}
			\left[D_{\operatorname{JS}}(Q)\right]^{(r)}.
		\end{aligned}
		\]
		Indeed, if \(s^{(r)}=0\), then
		$
		q_j^{(r)}=0
		$
		for every \(j\), and both the corresponding Euclidean term and
		$
		[D_{\operatorname{JS}}(Q)]^{(r)}
		$
		vanish. Hence
		\begin{equation}
			\label{eq:lipschitz-DJS-estimate}
			D_{\operatorname{JS}}(R)
			\preceqD
			\mathcal A_{A_{\operatorname{JS},s}}
			\bigl(D_{\operatorname{JS}}(Q)\bigr).
		\end{equation}
		
		Since \(\Neg_\Phi\) preserves hyperbolic mass, every iterate of
		\(P\in\PD(n)\) with \(s(P)=s\) remains in the same mass stratum.
		Therefore, iterating
		\eqref{eq:lipschitz-Dh-estimate},
		\eqref{eq:lipschitz-DE-estimate}, or
		\eqref{eq:lipschitz-DJS-estimate}, according to
		\(\nu\in\{h,E,\operatorname{JS}\}\), gives
		\[
		D_\nu\bigl(\Neg_\Phi^m(P)\bigr)
		\preceqD
		\mathcal A_{A_{\nu,s}}^{\,m}
		\bigl(D_\nu(P)\bigr).
		\]
		If \(\rho(A_{\nu,s})<1\), then, by
		Remark~\ref{rem:linear-comparison-map},
		\[
		\mathcal A_{A_{\nu,s}}^{\,m}(z)\longrightarrow0
		\qquad
		\text{for every }z\in\Dpos.
		\]
		It follows that
		\[
		D_\nu\bigl(\Neg_\Phi^m(P)\bigr)\longrightarrow0,
		\]
		which proves the corresponding assertion.
	\end{proof}
	
	\begin{theorem}
		\label{thm:practical-lipschitz-conditions}
		Let \(n\ge2\), let \(\Phi:\Dint{0}{1}\to\Dpos\) be an admissible
		hyperbolic generator, and fix \(s\in\{1,\eone,\etwo\}\).
		Suppose that
		\[
		\left|
		\frac{\partial\Phi^{(\ell)}}{\partial z^{(r)}}(z)
		\right|
		\le B_{\ell r},
		\qquad
		\Phi^{(\ell)}(q_j)\ge m_{\ell,s}>0
		\]
		for every active component \(\ell\), \(r\in\{1,2\}\),
		\(z\in\Dint{0}{1}\), and every 
		\(Q=(q_1,..., q_n)\in\PD(n)\) with \(s(Q)=s\).
		Then Proposition~\ref{prop:lipschitz-generated-negators-disequilibrium-decay}
		applies with $L_{\ell r}=B_{\ell r}$ and $M_{\ell,s}=n\,m_{\ell,s}$.
	\end{theorem}
	
	\begin{proof}
		The derivative bounds imply
		\eqref{eq:generator-lipschitz-condition} with
		\(L_{\ell r}=B_{\ell r}\). Moreover,
		\[
		Z_\ell^\Phi(Q)
		=
		\sum_{j=1}^{n}\Phi^{(\ell)}(q_j)
		\ge
		n\,m_{\ell,s},
		\]
		so \eqref{eq:normalizing-factor-lower-bound} holds with
		\(M_{\ell,s}=n\,m_{\ell,s}\).
		Therefore
		Proposition~\ref{prop:lipschitz-generated-negators-disequilibrium-decay}
		applies.
	\end{proof}
	
\begin{figure}[H]
		\centering
		\includegraphics[width=\textwidth]
		{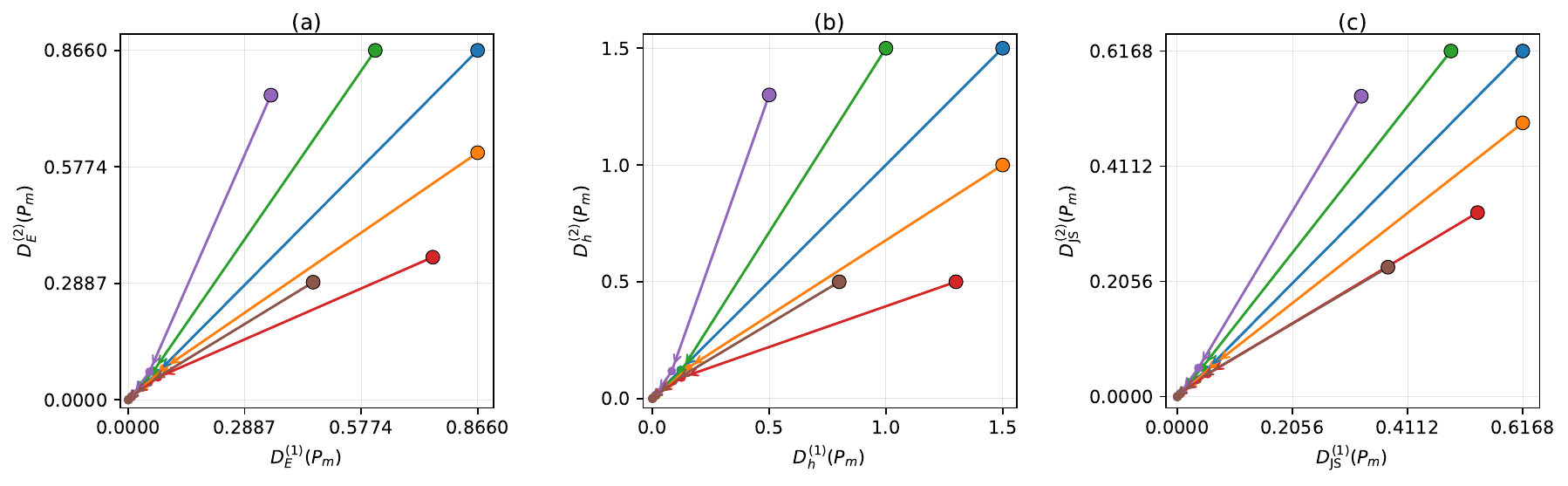}
		\caption{Decay of the hyperbolic-valued disequilibria associated with
			the iterates of \(\mathcal N_\gamma\) from
			Example~\ref{ex:lipschitz-matrix-exponential-negator}, with \(n=4\)
			and \(\gamma=0.5\). Panels \textup{(a)}, \textup{(b)}, and
			\textup{(c)} display the two idempotent components of
			\(D_E\), \(D_h\), and \(D_{\operatorname{JS}}\), respectively.}
		\label{fig:lipschitz-exponential-negator-disequilibrium-decay}
	\end{figure}

	\begin{example}
		\label{ex:lipschitz-matrix-exponential-negator}
		Let \(n\ge2\), and consider the stratum \(s(P)=1\). Let
		\(\gamma>0\), and define \(\Phi_\gamma:\Dint{0}{1}\to\Dpos\) by
		\[
		\Phi_\gamma(z)
		=
		\exp_{\mathbb D}
		\left(
		-\gamma
		\left[
		\left(
		\frac34 z^{(1)}+\frac14 z^{(2)}
		\right)\eone
		+
		\left(
		\frac13 z^{(1)}+\frac23 z^{(2)}
		\right)\etwo
		\right]
		\right).
		\]
		The generator \(\Phi_\gamma\) is order-reversing with respect to
		\(\preceqD\), since each idempotent component is non-increasing in
		both idempotent coordinates. Its components are strictly positive,
		and therefore its normalizing factors are positive. Hence
		Proposition~\ref{prop:generated-transformation-is-negator} implies that
		\(\mathcal N_\gamma=\Neg_{\Phi_\gamma}\) is a hyperbolic negator.
		
		On account of Theorem~\ref{thm:practical-lipschitz-conditions} set
		\[
		C=
		\begin{pmatrix}
			\frac34 & \frac14\\[1mm]
			\frac13 & \frac23
		\end{pmatrix},
		\]
		and let \(c_{\ell r}\) denote its \((\ell,r)\)-entry. For
		\(\ell,r\in\{1,2\}\), direct differentiation leads to
		\[
		\left|
		\frac{\partial \Phi_\gamma^{(\ell)}}{\partial z^{(r)}}(z)
		\right|
		=
		\gamma c_{\ell r}\Phi_\gamma^{(\ell)}(z)
		\le
		\gamma c_{\ell r},
		\]
		because
		$
		0<\Phi_\gamma^{(\ell)}(z)\le1
		$
		on \(\Dint{0}{1}\). Thus
		\eqref{eq:generator-lipschitz-condition} holds with
		$
		L_{\ell r}=\gamma c_{\ell r}
		$.
		
		Since the row sums of \(C\) are equal to \(1\), both linear
		combinations appearing in the exponent belong to \([0,1]\) for every
		\(z\in\Dint{0}{1}\). Therefore
		$
		\Phi_\gamma^{(\ell)}(z)\ge e^{-\gamma}
		$
		for \(\ell=1,2\). We may consequently take
		$
		m_{1,1}=m_{2,1}=e^{-\gamma}
		$,
		and Theorem~\ref{thm:practical-lipschitz-conditions} gives
		$
		M_{1,1}=M_{2,1}=ne^{-\gamma}
		$.
		
		The three matrices associated with the disequilibrium estimates are
		therefore
		\[
		A_{h,1}
		=
		\frac{2\gamma e^\gamma}{n}C,
		\qquad
		A_{E,1}
		=
		\frac{\gamma e^\gamma}{n}C,
		\qquad
		A_{\operatorname{JS},1}
		=
		\frac{\eta_n\gamma e^\gamma}{n}C.
		\]
		The eigenvalues of \(C\) are \(1\) and \(5/12\), so
		$
		\rho(C)=1
		$.
		Hence
		\[
		\rho(A_{h,1})
		=
		\frac{2\gamma e^\gamma}{n},
		\qquad
		\rho(A_{E,1})
		=
		\frac{\gamma e^\gamma}{n},
		\qquad
		\rho(A_{\operatorname{JS},1})
		=
		\frac{\eta_n\gamma e^\gamma}{n}.
		\]
		Thus, if
		\[
		\gamma e^\gamma
		<
		\frac{n}{\max\{2,\eta_n\}},
		\]
		all three spectral-radius conditions in
		Proposition~\ref{prop:lipschitz-generated-negators-disequilibrium-decay}
		hold simultaneously. Therefore, for every \(P\in\PD(n)\) with
		\(s(P)=1\),
		\[
		D_h\bigl(\mathcal N_\gamma^m(P)\bigr)\longrightarrow0,
		\qquad
		D_E\bigl(\mathcal N_\gamma^m(P)\bigr)\longrightarrow0,
		\qquad
		D_{\operatorname{JS}}
		\bigl(\mathcal N_\gamma^m(P)\bigr)\longrightarrow0.
		\]
		
		Figure~\ref{fig:lipschitz-exponential-negator-disequilibrium-decay}
		illustrates the decay of the three disequilibria for \(n=4\) and
		\(\gamma=0.5\).
	\end{example}
	
	\section*{Disclosure statement}
	
	The authors declare no competing interests.
	
	\section*{Funding}
	
	This research received no external or institutional funding.

\end{document}